\documentclass[english, 12pt,leqno]{amsart}

\usepackage{a4,amsfonts,amsmath,amssymb,amstext,epsfig,enumitem,color,url}

\newtheorem{thm}{Theorem}

\newtheorem{prop}[thm]{Proposition}

\theoremstyle{definition}

\newtheorem{rem}[thm]{Remark}
\newtheorem{defn}[thm]{Definition}

\newcommand{\Z}{\mathbf Z}

\newcommand{\R}{\mathbf R}

\newcommand{\Hi}{\mathcal H}

\newcommand{\Ki}{\mathcal K}

\newcommand{\id}{{\rm id}}

\title
{Spectra of infinite Cayley graphs,
\\
examples with pure band spectra}

\date{9 February 2024}

\author{Pierre de la Harpe}

\address{Pierre de la Harpe:
Section de math\'ematiques, 
Universit\'e de Gen\`eve,
C.P.~64, 
CH--1211 Gen\`eve 4.}
\email{Pierre.delaHarpe@unige.ch}



\thanks{The author acknowledges support of the Swiss NSF grant 200020-20040}

\begin{document}

\begin{abstract}
It is shown that there are groups $\Gamma$ with finite generating sets $S$
such that the adjacency operator of the Cayley graph ${\rm Cay}(\Gamma,S)$
is a disjoint union of $N$ intervals, for arbitrarily large integers $N$.
\end{abstract}

\maketitle

Let $G = (V, E)$ be a graph
with a countable (infinite or finite) set $V$ of vertices
and a set $E$ of edges which are unordered pairs of vertices;
graphs here are simple, i.e., without loops and multiple edges.
The degree of $v \in V$ is the number $d_v$ of edges of the form $\{u,v\}$
and we assume that $G$ is of bounded degree, that is $\max_{v \in V} d_v < \infty$.
The \textbf{adjacency operator} of $G$
is the linear operator $A_G$ acting on the Hilbert space $\ell^2(V)$ defined by
$$
(A_G \xi) (v) = \sum_{u \in V, \{u,v\} \in E} \xi(u)
\hskip.5cm \text{for} \hskip.5cm
\xi \in \ell^2(V) , \hskip.2cm v \in V .
$$
It is a bounded self-adjoint operator which appears for example in
\cite{CoSi--57, Bigg--93, Boll--98, BrHa--12} for finite graphs
and \cite{Kest--59, Moha--82, MoWo--89} for infinite graphs.
\par

Let $A$ be a bounded self-adjoint operator on a Hilbert space $\Hi$.
Denote by $\Sigma (A)$ the spectrum of $A$; 
it is a non-empty compact subset of the real line
 (non-empty because we assume $\Hi \ne \{ 0 \}$).
The operator $A$ has \textbf{pure band spectrum}
if its spectrum is the disjoint union of intervals,
$\Sigma (A) = \bigsqcup_{j=1}^N \mathopen[ a_j, b_j \mathclose]$
where $N \ge 1$ and $a_j < b_j < a_{j+1}$ for $j = 1, \hdots, N-1$ and $a_N < b_N$.
Graphs $G$ such that $A_G$ has pure band spectrum with $N \le 2$ are well-known,
for example Cayley graphs of finitely generated free abelian groups with $N=1$,
or free groups with $N = 1$ \cite{Kest--59},
% Grigorchuk--Zuk--2001 :
% The lamplighter group as a group generated by a 2-state automaton, and its spectrum
% Spectrum = interval (-1,1) but pure point spectrum.
or Grigorchuk groups with $N = 2$ \cite{DuGr--20}.
But, according to \cite{DuGr--20}, it has been an open question to find examples with $N \ge 3$.
% voir page 3 de DuGr--20.
(We repeat from \cite{DuGr--20} that it is another open question
whether there exists a Cayley graph such that the spectrum
of its adjacency operator is a Cantor set.)
The observation of this paper is that sums and products of graphs
provide a natural way to obtain graphs and \textbf{Cayley graphs}
with adjacency operators having pure band spectra.

\begin{thm}
% 1
\label{ThmMain}
For any $N \ge 1$, there exist Cayley graphs $G$ such that
the adjacency operator $A_G$ has a pure band spectrum,
consisting of $N$ disjoint intervals.
\end{thm}

See Section~\ref{SectionFinalEx} for precise examples,
which are Cayley graphs of groups which are products
$\Gamma_0 \times \Gamma_1 \times \cdots \times \Gamma_N$,
where $\Gamma_0$ is isomorphic to $\Z^d$ for some $d \ge 1$,
or to a free abelian group of rank $d \ge 2$,
and where $\Gamma_1, \hdots, \Gamma_N$ are finite groups.
\par

A graph $G = (V,E)$ has also a Laplacian $L_G$ and a Markov operator $M_G$
which act on $\ell^2(V)$. For regular graphs of some degree $k$, in particular for Cayley graphs,
we have $L_G = k - A_G$ and $M_G = 1 - \frac{1}{k} S_G$,
so that the spectra of $L_G$ and $M_G$ are images of $\Sigma(A_G)$
by affine transformations of $\R$.
The theorem above can therefore be repeated for $L_G$ and $M_G$.
\par

The proof of Theorem~\ref{ThmMain}, in Section~\ref{SectionFinalEx},
follows from standard results on spectra  of operators
acting on tensor products, see Section~\ref{SectionTPOp},
and on spectra on sums of graphs, see Section~\ref{SectionSPgraphs}.

\section{Spectrum of some operators acting on tensor products of Hilbert spaces}
% Section 1
\label{SectionTPOp}

Let $\Hi$ and $\Ki$ be Hilbert spaces,
$A$ a bounded self-adjoint operator on $\Hi$,
and $B$ a bounded self-adjoint operator on $\Ki$.
On the tensor product Hilbert space $\Hi \otimes \Ki$,
we consider the three operators
$$
A \otimes \id_\Ki + \id_\Hi \otimes B , \hskip.5cm
A \otimes B , \hskip.5cm
A \otimes \id_\Ki + \id_\Hi \otimes B + A \otimes B ,
$$
which are also bounded and self-adjoint.

\begin{prop}
% 2
\label{PropSpecTens}
Let $A$ and $B$ be as above.
\begin{enumerate}[label=(\roman*)]
%\addtocounter{enumi}{3}
\item\label{1DEPropSpecTens}
The spectrum of $A \otimes \id_\Ki + \id_\Hi \otimes B$ is
$$
\{ z \in \R \mid z = x+y
\hskip.2cm \text{for some} \hskip.2cm x \in \Sigma(A)
\hskip.2cm \text{and} \hskip.2cm y \in \Sigma(B) \} .
$$
\item\label{2DEPropSpecTens}
The spectrum of $A \otimes B$ is
$$
\{ z \in \R \mid z = x y
\hskip.2cm \text{for some} \hskip.2cm x \in \Sigma(A)
\hskip.2cm \text{and} \hskip.2cm y \in \Sigma(B) \} .
$$
\item\label{3DEPropSpecTens}
The spectrum of $A \otimes \id_\Ki + \id_\Hi \otimes B + A \otimes B$ is
$$
\{ z \in \R \mid z = x+y + xy
\hskip.2cm \text{for some} \hskip.2cm x \in \Sigma(A)
\hskip.2cm \text{and} \hskip.2cm y \in \Sigma(B) \} .
$$
\end{enumerate}
\end{prop}

There is one proof in \cite{Sche--69}, or in the older \cite{BrPe--66} for~\ref{2DEPropSpecTens}.
We sketch below another proof, which if needed could give more information
on spectral measures and multiplicities.
Before this, we recall the following facts.
A similar reminder can be found in \cite{Harp}.
\par

An operator (= bounded linear operator) $X$ on $\Hi$
and an operator $Y$ on $\Ki$ are \textbf{unitarily equivalent}
if there exists an operator $U \, \colon \Hi \to \Ki$ which is unitary (= a surjective isometry)
such that $Y = U X U^*$.
\par

Let $\Sigma$ be a non-empty metrizable compact space,
let $\mu$ be a finite measure on~$\Sigma$ with closed support equal to $\Sigma$,
let $\mathfrak m \, \colon \Sigma \to \{1, 2, 3, \hdots, \infty\}$
be a measurable function,
and let $\varphi \, \colon \Sigma \to \R \in L^\infty(\Sigma, \mu)$ be
an essentially bounded complex-valued function on $\Sigma$,
with essential supremum denoted by $\Vert \varphi \Vert_\infty$.
Denote by $\ell^2_\infty$ the Hilbert space
of square summable sequences $(z_j)_{j \ge 1}$ of complex numbers
and, for each $n \ge 1$, by $\ell^2_n$ the subspace
of sequences such that $z_j = 0$ for all $j \ge n+1$.
Let $L^2(\Sigma, \mu, \mathfrak m)$
be the Hilbert space of measurable functions $\xi \, \colon \Sigma \to \ell^2_\infty$
such that $\xi (x) \in \ell^2_{\mathfrak m (x)}$ for all $x \in \Sigma$
and $\int_\Sigma \Vert \xi (x) \Vert^2_{\ell^2_\infty} d\mu(x) < \infty$.
The \textbf{multiplication operator} $M_{\Sigma, \mu, \mathfrak m, \varphi}$
is the operator defined on the space $L^2(\Sigma, \mu, \mathfrak m)$ by
$$
(M_{\Sigma, \mu, \mathfrak m, \varphi} \xi) (x) = \varphi(x) \xi (x)
\hskip.5cm \text{for all} \hskip.2cm
\xi \in L^2(\Sigma, \mu, \mathfrak m)
\hskip.2cm \text{and} \hskip.2cm
x \in \Sigma .
$$
It is well-known that this is a bounded self-adjoint operator,
with norm $\Vert \varphi \Vert_\infty$ and with spectrum
the essential range of $\varphi$.
(The arguments to prove this are standard;
see for example Sections 4.20 to 4.28 in~\cite{Doug--72},
or any of \cite{AbKr--73, Abra--78, Krie--86}.)
A \textbf{straight multiplication operator} $M_{\Sigma, \mu, \mathfrak m}$
is an operator of this type in the particular case
of a compact subset $\Sigma$ of the real line
and of the function $\varphi$ given by the inclusion $\Sigma \subset \R$,
so that $(M_{\Sigma, \mu, \mathfrak m} ) (x) = x \xi (x)$
for all $\xi \in L^2(\Sigma, \mu, \mathfrak m)$
and $x \in \Sigma$.
\par

The Hahn--Hellinger multiplicity theorem establishes that:
\begin{enumerate}[label=(\Roman*)]
\item\label{1DEthmHH}
Two straight multiplication operators
$M_{\Sigma, \mu, \mathfrak m}$ and $M_{\Sigma', \mu' \mathfrak m'}$
are unitarily equivalent if and only if the three following conditions are satisfied:
\newline
$\Sigma' = \Sigma$, the measures $\mu'$ and $\mu$ are equivalent,
and $\mathfrak m'(x) = \mathfrak m(x)$ for $\mu$-almost all $x \in \Sigma$.
\item\label{2DEthmHH}
Any bounded self-adjoint operator $A$ on a separable Hilbert space $\Hi$
is unitarily equivalent to a straight multiplication operator $M_{\Sigma, \mu, \mathfrak m}$,
where $\Sigma$ is the spectrum of $A$.
\end{enumerate}
For a proof, we refer to 
\cite[Theorem 5.4.3]{Sim4--15} or to
\cite[Theorem 10.4.6]{BoSm--20} (which applies to unbounded self-adjoint operators).
The measure $\mu$ in~\ref{2DEthmHH} is called a scalar-valued spectral measure for $A$
(it is well-defined by $A$ up to equivalence of measures)
and the function $\mathfrak m$ is called the spectral multiplicity function of $A$
(is is well-defined up to equality $\mu$-almost everywhere).

\begin{proof}[A second proof of Proposition~\ref{PropSpecTens}]
Let $A$ and $B$ be as in the proposition.
Let $M_{\Sigma(A), \mu, \mathfrak m}$ and $M_{\Sigma(B), \nu, \mathfrak n}$
be two straight multiplication operators which are unitarily equivalent to $A$ and $B$ respectively.
Then $A \otimes \id_\Ki + \id_\Hi \otimes B$ is unitarily equivalent to
the multiplication operator (in general not a straight one)
$M_{\Sigma, \mu \times \nu, \mathfrak m \times \mathfrak n, \varphi_{i}}$
where $\Sigma$ is the product $\Sigma(A) \times \Sigma(B)$,
where $\mu \times \nu$ is the product measure on $\Sigma$,
where $(\mathfrak m \times \mathfrak n) (x,y) = \mathfrak m (x) \mathfrak n (y)$
and $\varphi_{i} (x,y) = x + y$ for all $(x,y) \in \Sigma$;
this operator acts on the Hilbert space
$L^2(\Sigma, \mu \times \nu, \mathfrak m \times \mathfrak n)
= L^2(\Sigma(A), \mu, \mathfrak m) \otimes L^2(\Sigma(B), \nu, \mathfrak n)$
of measurable functions $\zeta \, \colon \Sigma \to \ell^2_\infty \otimes \ell^2_\infty$
such that $\zeta(x,y) \in \ell^2_{\mathfrak m (x)} \otimes \ell^2_{\mathfrak n (y)}$
for all $(x,y) \in \Sigma$ and
$\int_\Sigma \Vert \zeta(x,y) \Vert^2 d(\mu \times \nu) (x,y) < \infty$.
Since the subset of $\R$ which appears in Claim~\ref{1DEPropSpecTens}
is the essential range of the function $\varphi_{i}$, 
this set is the spectrum of $M_{\Sigma, \mu \times \nu, \mathfrak m \times \mathfrak n, \varphi_{i}}$
by the reminder above,
and therefore also the spectrum of the unitarily equivalent operator
$A \otimes \id_\Ki + \id_\Hi \otimes B$.
This completes the proof of~\ref{1DEPropSpecTens}.
\par

The proofs of Claims~\ref{2DEPropSpecTens} and~\ref{3DEPropSpecTens}
are similar, using functions $\varphi_{ii}$ and $\varphi_{iii}$ defined by
$\varphi_{ii}(x,y) = xy$ and $\varphi_{iii}(x,y) = x + y + xy$.
\end{proof}

\section{Sums and products of graphs and Cayley graphs}
% Section 2
\label{SectionSPgraphs}

Let $G = (V, E)$ and $H = (W, F)$ be two simple graphs of bounded degree.
We recall the following definitions, as in~\cite[Section 2.5]{CvDS--80}.
% pages 65--66

\begin{defn}
% 3
\label{DefSPgraphs}
Let $G$ and $H$ be as above
\begin{enumerate}[label=(\roman*)]
\item\label{1DEDefSPgraphs}
The \textbf{sum} $G + H$ is the graph with vertex set $V \times W$ and edge set
$$
\begin{aligned}
E_{G + H} = & 
\Big\{ \{ (v,w), (v,w') \} \mid v \in V \hskip.2cm \text{and} \hskip.2cm \{w, w'\} \in F \Big \}
\\
& \bigcup 
\Big\{ \{ (v,w), (v',w) \} \mid \{v, v'\} \in E \hskip.2cm \text{and} \hskip.2cm w \in W \Big \}
\end{aligned}
$$
\item\label{2DEDefSPgraphs}
The \textbf{product} $G \times H$ is the graph with vertex set $V \times W$ and edge set
$$
E_{G \times H} = \Big\{ \{ (v,w), (v',w') \} \mid \{v, v'\} \in E \hskip.2cm \text{and} \hskip.2cm \{w, w'\} \in F \Big \}
$$
\item\label{3DEDefSPgraphs}
The \textbf{strong product} $G \times_s H$ is the graph with vertex set $V \times W$ and edge set
$$
E_{G \times_s H} = E_{G + H} \cup E_{G \times H} .
$$
\end{enumerate}
\end{defn}

Note that sums and strong products of connected graphs are connected,
but products need not be.
The following proposition is straightforward.

\begin{prop}
% 4
\label{PropSPgraphs}
Let $G = (V, E)$ and $H = (W, F)$ be two graphs as above,
with adjacency operators $A$ and $B$ respectively.
\begin{enumerate}[label=(\roman*)]
\item\label{1DEPropSPgraphs}
The adjacency operator of the sum $G + H$ is $A \otimes \id_{\ell^2(W)} + \id_{\ell^2(V)} \otimes B$.
\item\label{2DEPropSPgraphs}
The adjacency operator of the product $G \times H$ is $A \otimes B$.
\item\label{3DEPropSPgraphs}
The adjacency operator of the strong product $G \times_s H$
\newline
is $A \otimes \id_{\ell^2(W)} + \id_{\ell^2(V)} \otimes B + A \otimes B$.
\end{enumerate}
\end{prop}

From Propositions 2 and 4, we have the spectra of sum graphs, product graphs, and strong product graphs.
In the context of finite graphs, this is classical, see for example~\cite[Theorem 2.23]{CvDS--80},
% page 69
or~\cite[Section 1.4]{BrHa--12}.
\par

Here is an example of application.
Let $d \ge 1$ and let $Q_d$ the $1$-skeleton of the $d$-hypercube;
this graph is a Cayley graph of an elementary abelian $2$-group of order $2^d$,
and also the graph sum of $d$ copies
of the the graph $Q_1$ which has two vertices and one edge;
it follows from Propositions~\ref{PropSpecTens}~\ref{1DEPropSpecTens}
and~\ref{PropSPgraphs}~\ref{1DEPropSPgraphs}
that its spectrum has distinct eigenvalues
$\lambda_j = d-2j$, each of multiplicity $m_j = \binom{d}{j}$, for $j = 1, \hdots, d$;
this is of course well-known \cite[Result 21a]{Bigg--93}.

\begin{rem}
% 5
\label{RemSumsCayley}
For two Cayley graph $G = {\rm Cay}(\Gamma,S)$ and $H = {\rm Cay}(\Delta,T)$,
the sum $G+H$ is the Cayley graph of the group $G \times H$
with respect to the generating set  $S \cup T \Doteq ( S \times \{e_H\} ) \cup ( \{e_G\} \times T )$
and the strong product $G \times_s H$ is the Cayley graph of $G \times H$
with respect to the generating set $S \cup T \cup (S \times T)$.
Note that $S \times T$ need not generate $\Gamma \times \Delta$
(example: $\Gamma = \Delta = \Z / 2\Z$ andd $\vert S \vert = \vert T \vert = 1$);
when it does, $G \times_s H$ is the Cayley graph of $\Gamma \times \Delta$
with respect to $S \times T$.
\end{rem}

\section{Examples and a proof of Theorem~1}
% Section 3
\label{SectionFinalEx}

Let $n$ be a positive integer, $K_n$ the complete graph on $n$ vertices,
and $A_{(n)}$ the adjacency operator of $K_n$.
The spectrum of $A_{(n)}$ consists of the simple eigenvalue $n-1$
and of the eigenvalue $-1$ of multiplicity $n-1$.
The graph $K_n$ is the Cayley graph of a group $G$ of order $n$
generated by the set $G \smallsetminus \{e_G\}$.
\par

Let $N$ be a positive integer and let $NK_n$
be the sum of $N$ copies of $K_n$.
It follows from Proposition~\ref{PropSPgraphs} that
the distinct eigenvalues of the adjacency operator of $NK_n$
are $-(N-j) + j(n-1)$ for $j = 0, 1, \hdots, N$
(we don't need their multiplicities here).
Note that the gap between two consecutive eigenvalues is $n$.
\par

Let $L$ be the infinite line, with vertex set $\Z$
and edge set $\Big\{ \{ n, n+1\} \mid n \in \Z \Big\}$.
Using Fourier transform, it is standard to compute the spectrum
of the adjacency operator of $L$, which is the interval $\mathopen[ -2, 2 \mathclose]$;
the computation is repeated in \cite{Harp}.
By Remark~\ref{RemSumsCayley}, the sum $L + NK_n$
is a Cayley graph of a group product of $\Z$ and $N$ groups of order $n$.
It is now an immediate consequence of
Propositions~\ref{PropSPgraphs} \ref{1DEPropSPgraphs} 
and~\ref{PropSpecTens} \ref{1DEPropSpecTens} that we have:

\begin{prop}
% 6
\label{PropFin}
Let $n \ge 5$ and $N \ge 1$.
The graph $L + NK_n$
is a Cayley graph of which the spectrum of the adjacency operator
is the disjoint union of $N+1$ intervals of length $4$ centered at the points
$-(N-j) + j(n-1)$, where $j = 0, 1, \hdots, N$.
\end{prop}

We can replace $L$ by any Cayley graph with spectrum an interval,
as long as $n$ is large enough, for example by the standard Cayley graph
of the free abelian group of rank $d$,
which has a spectrum $\mathopen[ -2d, 2d \mathclose]$ (where $d \ge 1$),
or by a free group of rank $d$,
which has a spectrum $\mathopen[ -2\sqrt{d-1}, 2 \sqrt{d-1} \mathclose]$
(where $d \ge 2$).
\par

We can also replace $K_n$ by the complete bipartite graph $K_{n,n}$
which is the Cayley graph of a group $G$ or order $2n$
which has a subgroup $H$ of order $n$,
with respect to the generating set $G \smallsetminus H$.
The eigenvalues of $K_{n,n}$ are the simple eigenvalues $-n$ and $n$,
and $0$ of multipicity $2n-2$.
For $n \ge 5$, the spectrum of $L + NK_{n,n}$
is the disjoint union of 2N+1 intervals of length $4$ centered
at the $2N+1$  points $jN$, where $-N \le j \le N$.
\par

Note that all groups which can be used for Proposition~\ref{PropFin} have torsion.
Let $\Gamma$ be a countable torsion-free group.
If the reduced C$^*$-algebra $C^*_r(\Gamma)$ has no non-trivial projection,
the spectrum of a Cayley graph of $\Gamma$ is connected,
and thus cannot be a spectrum with $N \ge 2$ bands, in contrast with the proposition.
The Kadison--Kaplansky conjecture is the statement that
$C^*_r(\Gamma)$ has no non-trivial projection for any countable torsion-free group;
it is known to  hold for large classes of groups,
including torsion-free hyperbolic groups \cite{Pusc--02},
their subgroups \cite{ MiYu--02},
and torsion-free amenable groups
(as a consequence of heavy results on the Baum--Connes conjecture~\cite[Corollary~9.2]{HiKa--01}).


\begin{thebibliography}{AkGW--08}

\bibitem[Abra--78]{Abra--78}
M.B.\ Abrahamse,
\emph{Multiplication operators}.
In ``Hilbert space operators (Proc.\ Conf., Calif.\ State Univ., Long Beach, Calif., 1977)'',
Lectures Notes in Math.\ \textbf{693}, Springer (1978), 17--36.

\bibitem[AbKr--73]{AbKr--73}
M.B.\ Abrahamse and T.L.\ Kriete,
\emph{The spectral multiplicity of a multiplication operator}.
Indiana U.\ Math.\ Jour.\ \textbf{22} (1973), 845--857.

\bibitem[Bigg--93]{Bigg--93}
N.\ Biggs,
\emph{Algebraic graph theory}, Second edition.
Cambridge Univ.\ Press, 1993.

\bibitem[Boll--98]{Boll--98}
B.\ Bollob\'as,
\emph{Modern graph theory}.
Graduate Texts in Math.\ \textbf{198},
Springer 1998.

\bibitem[BoSm--20]{BoSm--20}
V.I.\ Bogachev and O.G.\ Smolyanov,
\emph{Real and functional analysis}.
Moscow Lectures \textbf{4}, Springer, 2020.

\bibitem[BrHa--12]{BrHa--12}
A.E.\ Brouwer and W.H.\ Haemers,
\emph{Spectra of graphs}.
Universitext, Springer, 2012.

\bibitem[BrPe--66]{BrPe--66}
A.\ Brown and C.\ Pearcy,
\emph{Spectra of tensor products of operators}.
Proc.\ Amer.\ Math.\ Soc.\ \textbf{17} (1966), no.\ 1, 162--166.

\bibitem[CoSi--57]{CoSi--57}
L.\ Collatz and U.\ Sinogowitz,
\emph{Spektren endlicher Grafen}.
Abh.\ Math.\ Sem.\ Univ.\ Hamburg \textbf{21} (1957), 63--77.

\bibitem[CvDS--80]{CvDS--80}
D.M.\ Cvetkovi\'c, M.\ Doob, and H.\ Sachs,
\emph{Spectra of graphs. Theory and applications.}
Academic Press, 1980.

\bibitem[Doug--72]{Doug--72}
R.G.\ Douglas,
\emph{Banach algebra techniques in operator theory}.
Academic Press, 1972.
[Second Edition, Springer, 1998.]

\bibitem[DuGr--20]{DuGr--20}
A.\ Dudko and R.\ Grigorchuk,
\emph{On the question ''can one hear the shape of a group''
and a Hulanicki type theorem for graphs}.
Israel J.\ Math. \textbf{237} (2020), 53--74.

\bibitem[Harp]{Harp}
P.\ de la Harpe,
\emph{Spectral multiplicity functions of adjacency operators of graphs
and cospectral infinite graphs}.
ArXiv:2308.04339v1, 8 Aug 2023.

\bibitem[HiKa--01]{HiKa--01}
N.\ Higson and G.\ Kasparov,
\emph{E-theory and KK-theory for groups which act
properly and isometrically on Hilbert space}.
Invent.\ Math.\ \textbf{144} (2001), 23--74.

\bibitem[Kest--59]{Kest--59}
H.\ Kesten,
\emph{Symmetric random walks on groups}.
Trans.\ Amer.\ Math.\ Soc.\ \textbf{92} (1959), no.\ 2, 336--354.

\bibitem[Krie--86]{Krie--86}
T.L.\ Kriete,
\emph{An elementary approach to the multiplicity theory of multiplication operators}.
Rocky Mountain J.\ Math.\ \textbf{16} (1986), no.\ 1, 23--32.

\bibitem[MiYu--02]{MiYu--02}
I.\ Mineyev and G.\ Yu,
\emph{The Baum--Connes conjecture for hyperbolic groups}.
Invent.\ Math.\ \textbf{149} (2002), no.\ 1, 97--122.

\bibitem[Moha--82]{Moha--82}
B.\ Mohar,
\emph{The spectrum of an infinite graph}.
Linear Algebra Appl.\ \textbf{48} (1982), 245--256.

\bibitem[MoWo--89]{MoWo--89}
B.\ Mohar and W.\ Woess,
\emph{A survey on spectra of infinite graphs}.
Bull.\ London Math.\ Soc.\ \textbf{21} (1989), no\ 3, 209--234.

\bibitem[Pusc--02]{Pusc--02}
M.\ Puschnigg,
\emph{The Kadison--Kaplansky conjecture for word-hyperbolic groups}.
Invent.\ Math.\ \textbf{149} (2002), 153--194.

\bibitem[Sche--69]{Sche--69}
M.\ Schechter,
\emph{On the spectra of operators on tensor products}.
J.\ Functional Analysis \textbf{4} (1969), 95--99.

\bibitem[Sim4--15]{Sim4--15}
B.\ Simon,
\textit{Operator theory},
Compr.\ Course Anal., Part 4,
Amer.\ Math.\ Soc., 2015.

\end{thebibliography}
\end{document}